\documentclass[a4paper]{article}
\usepackage{amsmath,amssymb,mathrsfs,graphicx,psfrag}
\allowdisplaybreaks
\DeclareMathOperator{\re}{Re}
\DeclareMathOperator{\im}{Im}
\providecommand{\keywords}[1]{Keywords: #1}
\providecommand{\MSC}[1]{MSC: #1}
\begin{document}
\title{Method of fundamental solutions for the problem of doubly-periodic potential flow}
\author{Hidenori Ogata\footnote{Department of Computer and Network Engineering, Graduate School of Informatics and Engineering, The University of Electro-Communications, 1-5-1 Chofugaoka, Chofu, Tokyo, 182-8585, Japan. (e-mail) {\tt ogata@im.uec.ac.jp}}}
\maketitle
\begin{abstract}
 In this paper, we propose a method of fundamental solutions for the problem of 
 two-dimensional potential flow in a doubly-periodic domain. 
 The solution involves a doubly-periodic function, to which it is difficult to give an approximation 
 by the conventional method of fundamental solutions. 
 We propose to approximate it by a linear combination of the periodic fundamental solutions, 
 that is, complex logarithmic potentials with sources in a doubly-periodic array constructed using the theta functions. 
 Numerical examples show the effectiveness of our method. 
\end{abstract}
\keywords{method of fundamental solutions, potential problem, double periodicity, periodic fundamental solution, theta function, elliptic function}

\medskip

\noindent
\MSC{65N80,65E05}
\section{Introduction}
\label{sec:introduction}
The method of fundamental solutions, or the charge simulation method, 
\cite{FairweatherKarageorghis1998,Murashima1983} 
is a fast numerical solver for potential problems
\begin{equation}
 \label{eq:potential-problem}
  \begin{cases}
   - \triangle u = 0 & \mbox{in} \ \ \mathscr{D}
   \\
   u = f & \mbox{on} \ \ \partial\mathscr{D}, 
  \end{cases}
\end{equation}
where $\mathscr{D}$ is a domain in the $n$-dimensional Euclidean space $\mathbb{R}^n$, 
and $f$ is a function given on $\partial\mathscr{D}$. 
In two dimensional problems $(n=2)$, equalizing the Euclid plane $\mathbb{R}^2$ 
with the complex plane $\mathbb{C}$, 
the method of fundamental solutions gives an approximate solution of 
(\ref{eq:potential-problem}) of the form 
\footnote{The scheme shown here is the invariant scheme \cite{Murota1993,Murota1995} proposed by Murota.}
\begin{equation}
 \label{eq:mfs-solution0}
  u(z) \simeq u_N(z) = 
  Q_0 - \frac{1}{2\pi}\sum_{j=1}^{N}Q_j\log|z - \zeta_j|
  \quad ( \: z = x + \mathrm{i}y \: ), 
\end{equation}
where $Q_0, Q_1, \ldots, Q_N$ are unknown real coefficients such that 
\begin{equation}
 \label{eq:q-sum-zero0}
  \sum_{j=1}^{N}Q_j = 0, 
\end{equation}
and $\zeta_1, \ldots, \zeta_N$ are points given in $\mathbb{C}\setminus\overline{\mathscr{D}}$. 
We call $Q_j$ the ^^ ^^ charges'' and $\zeta_j$ the ^^ ^^ charge points''. 
We remark that the approximate solution $u_N(z)$ exactly satisfies the Laplace equation in $\mathscr{D}$. 
Regarding the boundary condition, we pose the collocation condition on $u_N(z)$, namely, 
we assume that $u_N(z)$ satisfies the equations
\begin{equation}
 \label{eq:collocation-cond0}
  u_N(z_i) = f(z_i), \quad i = 1, \ldots, N 
\end{equation}
with $z_1, \ldots, z_N$ given on the boundary $\partial\mathscr{D}$, 
which are called the ^^ ^^ collocation points''. 
The equations (\ref{eq:collocation-cond0}) are rewritten as
\begin{equation}
 \label{eq:collocation-cond02}
  Q_0 - \frac{1}{2\pi}\sum_{j=1}^{N}Q_j\log|z_i - \zeta_j| = f(z_i), \quad 
  i = 1, \ldots, N, 
\end{equation}
which, together with (\ref{eq:q-sum-zero0}), form a system of linear equations 
with respect to $Q_0, Q_1, \ldots, Q_N$. 
We determine the unknowns $Q_j$ by solving the linear system (\ref{eq:q-sum-zero0}) and 
(\ref{eq:collocation-cond02}) and obtain the approximate solution $u_N(z)$. 
The method of fundamental solutions has the advantages that it is easy to program, 
its computational cost is low, and it shows fast convergence such as exponential convergence 
\cite{KatsuradaOkamoto1988} under some condition. 
It was first used for studies of electric field problems 
\cite{SingerSteinbiglerWeiss1974,Steinbigler-dissertation1969},  
and now it is widely used in science and engineering, 
for example, problem of scattering of earthquake waves \cite{Sanchez-SezmaRosenblueth1979}.

The method of fundamental solutions is also used for the approximation of complex analytic functions. 
Let $f(z)$ be an analytic function in a domain $\mathscr{D}\subset\mathbb{C}$. 
The real part of $f(z)$, which is a harmonic function in $\mathscr{D}$, can be approximated 
using the form (\ref{eq:mfs-solution0}), and the imaginary part of $f(z)$, which is the conjugate harmonic function 
of $\re f(z)$, is approximated using 
\begin{equation*}
 - \frac{1}{2\pi}\sum_{j=1}^{N}Q_j\arg(z - \zeta_j). 
\end{equation*}
Then, the analytic function $f(z)$ is approximated using a linear combination of the complex logarithmic functions 
\begin{equation}
 \label{eq:complex-mfs-solution}
  Q_0 - \frac{1}{2\pi}\sum_{j=1}^{N}Q_j\log(z - \zeta_j). 
\end{equation}
From this point of view, Amano applied the method of fundamental solution to numerical conformal mappings 
\cite{Amano1994,Amano1998,AmanoOkanoOgataSugihara2012}. 

\begin{figure}[htbp]
 \begin{center}
  \psfrag{w}{$\omega_1$}
  \psfrag{W}{$\omega_2$}
  \psfrag{0}{$D_{00}$}
  \psfrag{1}{$D_{10}$}
  \psfrag{2}{$D_{01}$}
  \psfrag{3}{$D_{11}$}
  \psfrag{D}{$\mathscr{D}$}
  \includegraphics[width=0.5\textwidth]{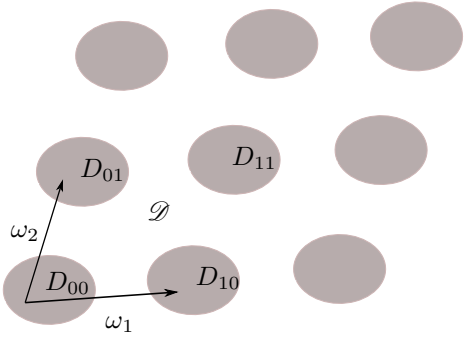}
 \end{center}
 \caption{Flow domain $\mathscr{D}$ with an infinite doubly-periodic array of obstacles 
 $D_{mn}, m,n\in\mathbb{Z}$.}
 \label{fig:periodic-domain}
\end{figure}

In this paper, we examine the problem of two-dimensional potential flow 
past an infinite doubly-periodic array of obstacles as shown in Figure \ref{fig:periodic-domain}. 
A two-dimensional potential flow is characterized by a complex velocity potential $f(z)$ 
\cite{Milne-Thomson2011}, 
which is an analytic function in the flow domain $\mathscr{D}$ such that it gives the velocity field 
$\boldsymbol{v}=(u,v)$ by $f^{\prime}(z) = u - \mathrm{i}v$, and its imaginary part satisfies 
the boundary condition 
\begin{equation}
 \label{eq:boundary-cond0}
  \im f = \mbox{constant} \quad \mbox{on} \ \ \partial\mathscr{D}. 
\end{equation}
Physically, the boundary condition (\ref{eq:boundary-cond0}) means that the fluid flows along the boundary 
$\partial\mathscr{D}$ since the contour lines of $\im f$ are the streamlines. 
Therefore, in order to obtain the potential flow of our problem, 
we have to find an analytic function $f(z)$ in $\mathscr{D}$ satisfying the boundary condition 
(\ref{eq:boundary-cond0}). 

We have, however, one problem in applying the method of fundamental solutions to our problem. 
The velocity field $\boldsymbol{v}$ is obviously a doubly-periodic function 
due to the the double periodicity of the flow domain $\mathscr{D}$. 
Then, the complex velocity potential $f(z)$ involves a doubly-periodic function, and 
it is difficult to approximate it using the form (\ref{eq:complex-mfs-solution}) 
of the conventional method. 
To overcome this challenge, we propose to solve our problem 
using a doubly-periodic fundamental solution, that is, a logarithmic potential with a doubly-periodic array 
of sources. 
We construct a doubly-periodic fundamental solution using the theta functions 
and approximate the complex velocity potential by a linear combination of the doubly-periodic fundamental solutions. 

The previous works related to this paper are as follows. 
As studies of problems of periodic flow, Zick and Homsy \cite{ZickHomsy1982} proposed an integral equation method for 
three-dimensional Stokes flow problems with a three-dimensional periodic array of spheres, 
where the solution is given by an integral including the periodic fundamental solution. 
Greengard and Kropinski \cite{GreengardKropinski2004} proposed an integral equation method for 
two-dimensional Stokes flow problems in doubly-periodic domains, where an approximate solution is given 
as a complex variable formulation and the fast multipole method is used. 
Liron \cite{Liron1978} studied Stokes flow due to infinite array of Stokeslets and 
applied it to the problems of fluid transport by cilia.
As studies on methods of fundamental solutions applied to problems with periodicity, 
Ogata et al. proposed a method of numerical conformal mappings of complex domains with a single periodicity 
\cite{OgataOkanoAmano2002}, 
where the mapping function, an analytic function involving a singly-periodic function, is approximated by 
the method of fundamental solution using the singly-periodic fundamental solutions.  
They proposed a method of fundamental solutions also to the problems of two-dimensional periodic Stokes flows 
\cite{OgataAmanoSugiharaOkano2003,OgataAmano2010}, where the solutions are approximated by the periodic 
fundamental solutions of the Stokes equation, that is, the Stokes flows induced by a periodic array of point forces. 
The author proposed a method of fundamental solutions for the two-dimensional elasticity problem with one-dimensional 
periodicity \cite{Ogata2008}, where the solution is approximated using the periodic fundamental solutions 
of the elastostatic equation, that is, the displacements induced by concentrated forces in a periodic array. 

The remainder of this paper is structured as follows. 
Section \ref{sec:method-fundamental-solutions} proposes a method of fundamental solutions for our problem of 
potential flow with double periodicity. 
Section \ref{sec:numerical-example} shows some numerical examples which show the effectiveness of our method. 
In Section \ref{sec:concluding-remark}, we conclude this paper and present problems related to future studies. 
\section{Method of fundamental solutions}
\label{sec:method-fundamental-solutions}
\begin{figure}[htbp]
 \begin{center}
  \psfrag{D}{$\mathscr{D}_0$}
  \psfrag{z}{$z_0$}
  \psfrag{w}{$\omega_1$}
  \psfrag{W}{$\omega_2$}
  \includegraphics[width=0.5\textwidth]{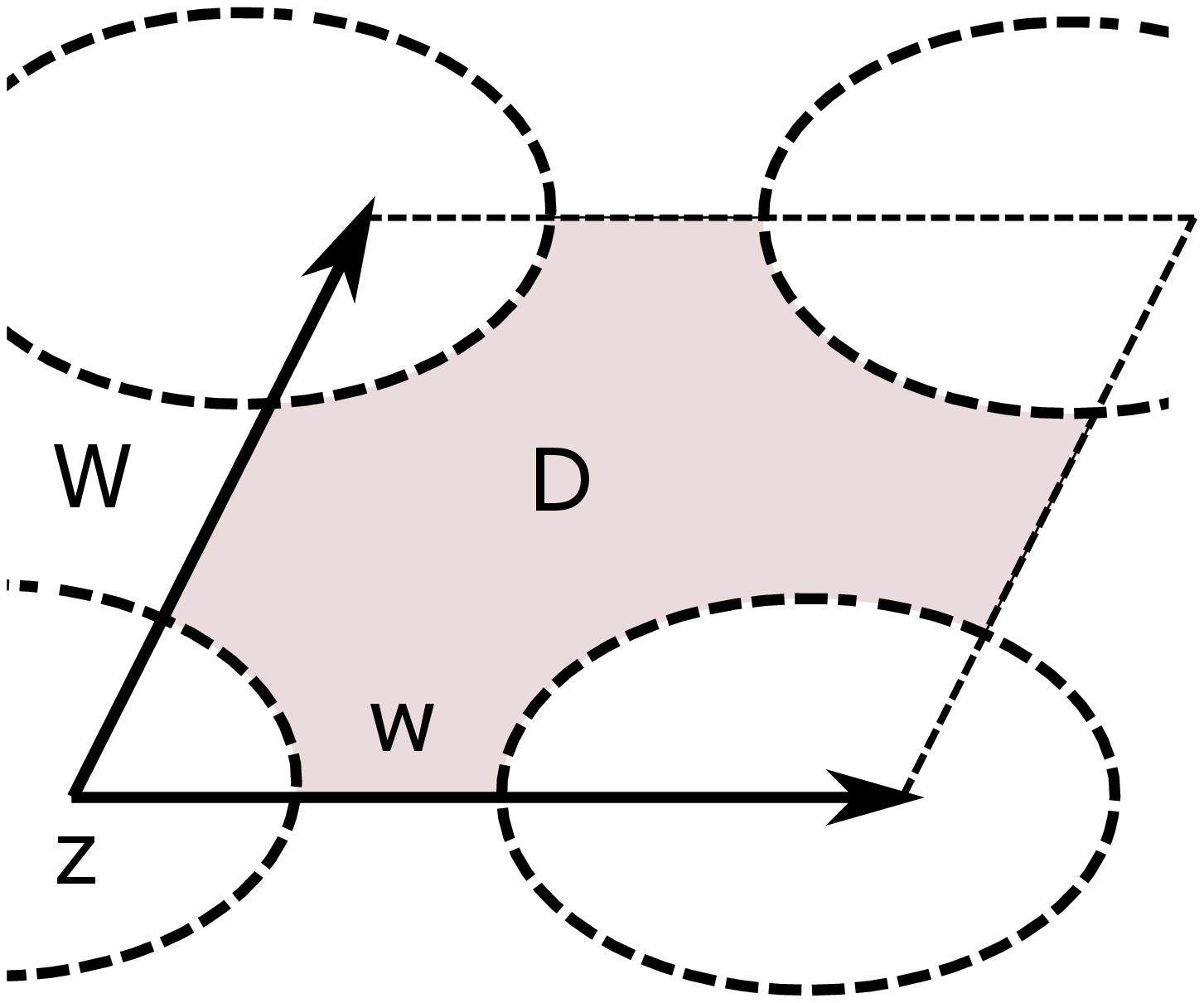}
 \end{center}
 \caption{The set $\mathscr{D}_0$.}
 \label{fig:set-D0}
\end{figure}

We consider a potential flow past the doubly-periodic array of obstacle. 
The flow domain is mathematically given by 
\begin{equation*}
 \mathscr{D} = 
  \mathbb{C}\setminus\bigcup_{m,n\in\mathbb{Z}}\overline{D_{mn}}, 
\end{equation*}
where $D_{0}$ is one of the obstacles, a simply-connected domain in $\mathbb{C}$ and 
\begin{equation*}
 D_{mn} = 
  \left\{ \: z + m\omega_1 + n\omega_2 \: \left| \: z \in D_{00} \: \right.\right\}, 
  \quad 
  m, n \in \mathbb{Z}
\end{equation*}
with complex numbers $\omega_1$ and $\omega_2$ giving the periods of the obstacle array such that 
$\im(\omega_2/\omega_1)>0$ and 
$\overline{D_{mn}}\cap\overline{D_{m^{\prime}n^{\prime}}}=\emptyset$ if 
$(m,n)\neq(m^{\prime},n^{\prime})$. 
We assume that the spatial average of the velocity field $\boldsymbol{v}=(u,v)$ is 
a uniform flow, in whose direction the real axis is taken, that is, 
\begin{equation}
 \label{eq:uniform-flow}
  \langle\boldsymbol{v}\rangle = 
  \frac{1}{|\mathscr{D}_0|}\iint_{\mathscr{D}_0}\boldsymbol{v}\mathrm{d}x\mathrm{d}y = (U, 0),
\end{equation} 
where $U$ is the magnitude of the velocity of the unit flow, 
$\mathscr{D}_0$ is the set defined by 
\begin{align}
 \nonumber  
 \mathscr{D}_0 = \: & 
 \mbox{(a period parallelogram)}\setminus\bigcup_{m,n\in\mathbb{Z}}\overline{D_{mn}}
 \\
 = \: & 
 \left\{ \: 
 z_0 + a_1\omega_1 + a_2\omega_2 \: \left|\right.\: 0\leq a_1, a_2\leq 1 \: 
 \right\} 
 \setminus\bigcup_{m,n\in\{0,1\}}\overline{D_{mn}}
\end{align}
with a point
\begin{math}
 z_0 \in D_{00}, 
\end{math}
(see Figure \ref{fig:set-D0}), and $|\mathscr{D}_0|$ is the area of $\mathscr{D}_{0}$. 

As mentioned in the previous section, 
a two-dimensional potential flow is characterized by a complex velocity potential, 
$f(z)$, which is an analytic function in the flow domain $\mathscr{D}$ 
such that it gives the velocity field $\boldsymbol{v}=(u,v)$ by 
\begin{equation}
 \label{eq:complex-velocity}
  f^{\prime}(z) = u - \mathrm{i}v,
\end{equation}
and its imaginary part is constant on the boundary of the domain $\mathscr{D}$. 
The condition (\ref{eq:uniform-flow}) is rewritten as 
\begin{equation}
 \label{eq:uniform-flow2}
  \frac{1}{|\mathscr{D}_0|}\iint_{\mathscr{D}_0}f^{\prime}(z)\mathrm{d}x\mathrm{d}y = U. 
\end{equation}
Therefore, our potential flow problem is mathematically reduced to the problem 
to find an analytic function $f(z)$ in $\mathscr{D}$ satisfying 
\begin{equation}
 \label{eq:boundary-condition}
  \im f = \mbox{constant} \quad \mbox{on} \ \ \partial\mathscr{D}
\end{equation}
and (\ref{eq:uniform-flow2}). 

In solving our problem by the method of fundamental solution method, however, 
we have one problem. 
The complex velocity potential $f(z)$ involves a periodic function since the velocity field 
$\boldsymbol{v}=(u,v)$ given by (\ref{eq:complex-velocity}) is obviously a periodic function, 
and it is impossible to approximate $f(z)$
by the conventional method using the form (\ref{eq:mfs-solution0}). 
To overcome this challenge, we propose to approximate the complex velocity potential $f$ by 
\begin{equation}
 \label{eq:mfs-solution}
  f(z) \simeq f_N(z) = 
  Uz - 
  \frac{\mathrm{i}}{2\pi}\sum_{j=1}^{N}Q_j
  \left\{
  \log\vartheta_1\left(\left.\frac{z-\zeta_j}{\omega_1}\right|\tau\right)
  - u_j z
  \right\}, 
\end{equation}
where $\vartheta_1(v|\tau)$ is the theta function \cite{Armitage-Eberlein2006} 
\footnote{
In \cite{Armitage-Eberlein2006}, the theta function $\vartheta_1(z|\tau)$ is defined by 
\begin{equation*}
 \begin{aligned}
  \vartheta_1(z|\tau) = \: & 
  2\sum_{n=0}^{\infty}q^{(n+1/2)^2}\sin(2n+1)z 
  \\ 
  = \: & 
  2q^{1/4}\sin z\prod_{n=1}^{\infty}
  (1 - q^{2n})
  (1 - 2q^{2n}\cos 2z + q^{4n}), 
  \quad q = \mathrm{e}^{\mathrm{i}\pi\tau}.
 \end{aligned}
\end{equation*}
}
\begin{equation*}
 \begin{aligned}
  \vartheta_1(v|\tau) = \: & 
  2\sum_{n=0}^{\infty}q^{(n+1/2)^2}\sin(2n+1)\pi v
  \\
  = \: & 
  2q^{1/4}\sin\pi v \prod_{n=1}^{\infty}
  (1 - q^{2n})(1 - 2q^{2n}\cos 2\pi v + q^{4n})
 \end{aligned}
\end{equation*}
with $\tau = \omega_2/\omega_1$ and $q = \mathrm{e}^{\mathrm{i}\pi\tau}$, 
\begin{equation}
 \label{eq:uj}
  u_j = 
  \frac{1}{|\mathscr{D}_0|}\iint_{\mathscr{D}_0}
  \frac{1}{\omega_1}
  \frac{\vartheta_1^{\prime}((z-\zeta_j)/\omega_1|\tau)}{\vartheta_1((z-\zeta_j)/\omega_1|\tau)}
  \mathrm{d}x\mathrm{d}y, 
  \quad 
  j = 1, \ldots, N, 
\end{equation}
$\zeta_1, \ldots, \zeta_N$ are points given in $D_{00}$ and 
$Q_1, \ldots, Q_N$ are unknown real coefficients such that 
\begin{equation}
 \label{eq:q-sum-zero}
  \sum_{j=1}^{N}Q_j = 0.
\end{equation}
We call $Q_1, \ldots, Q_N$ the ^^ ^^ charge'' and 
$\zeta_1, \ldots, \zeta_N$ the ^^ ^^ charge points''.
The theta function $\vartheta_1(v|\tau)$ is an entire function which has simple zeros 
at $m+n\tau$, $m,n\in\mathbb{Z}$ and satisfies the pseudo-periodicity
\begin{equation}
 \label{eq:theta1-pseudo-periodicity}
  \vartheta_1(v+1|\tau) = - \vartheta_1(v|\tau), \quad 
  \vartheta_1(v+\tau|\tau) = - q^{-1}\mathrm{e}^{-2\pi\mathrm{i}v}\vartheta_1(v|\tau). 
\end{equation}
Therefore, the functions $-(1/(2\pi))\log\vartheta_1((z-\zeta_j)/\omega_1|\tau)$ 
appearing on the rightmost side of (\ref{eq:mfs-solution}) can be regarded 
as the complex logarithmic potential with sources at the points 
\begin{math}
 m\omega_1 + n\omega_2 + \zeta_j, \ m, n \in \mathbb{Z},  
\end{math}
and it is a periodic fundamental solution of the two-dimensional Poisson equation 
\footnote{
Hasimoto \cite{Hasimoto2008} presented a periodic fundamental solution of the two-dimensional 
Poisson equation in terms of the Weierstrass elliptic functions
}. 
The term $-u_j z$ is added in each term of the rightmost side of (\ref{eq:mfs-solution}) 
so that $f_N(z)$ satisfies the condition (\ref{eq:uniform-flow2}). 

The approximate potential $f_N(z)$ satisfies the pseudo-periodicity
\begin{align}
 \label{eq:pseudo-periodicity-potential1}
 f_N(z+\omega_1)-f_N(z) = \: & 
 \omega_1\left(U + \frac{\mathrm{i}}{2\pi}\sum_{j=1}^{N}Q_ju_j\right), 
 \\
 \label{eq:pseudo-periodicity-potential2}
 f_N(z+\omega_2)-f_N(z) = \: & 
 U\omega_2 + 
 \sum_{j=1}^{N}Q_j
 \left( \frac{\zeta_j}{\omega_1} + \frac{\mathrm{i}}{2\pi}\omega_2 u_j \right). 
\end{align}
In fact, 
\begin{align*}
 & 
 f_N(z+\omega_1) - f_N(z) 
 \\
 = \: & 
 U\omega_1 - 
 \frac{\mathrm{i}}{2\pi}\sum_{j=1}^{N}Q_j
 \left\{
 \log\vartheta_1\left(\left.\frac{z-\zeta_j}{\omega_1}+1\right|\tau\right) - 
 \log\vartheta_1\left(\left.\frac{z-\zeta_j}{\omega_1}\right|\tau\right) - u_j\omega_1
 \right\}
 \\
 = \: & 
 U\omega_1 - 
 \frac{\mathrm{i}}{2\pi}\sum_{j=1}^{N}Q_j\left\{\log(-1)-u_j\omega_1\right\}
 = 
 \omega_1\left(U+\frac{\mathrm{i}}{2\pi}\sum_{j=1}^{N}Q_ju_j\right), 
 \intertext{where (\ref{eq:theta1-pseudo-periodicity}) is used on the second equality and 
(\ref{eq:q-sum-zero}) on the third equality, and }
 & 
 f_N(z+\omega_2) - f_N(z)
 \\
 = \: & 
 U\omega_2 - 
 \frac{\mathrm{i}}{2\pi}\sum_{j=1}^{N}Q_j
 \left\{
 \log\vartheta_1\left(\left.\frac{z-\zeta_j}{\omega_1}+\tau\right|\tau\right)
 - 
 \log\vartheta_1\left(\left.\frac{z-\zeta_j}{\omega_1}\right|\tau\right)
 - 
 u_j\omega_2
 \right\}
 \\
 = \: & 
 U\omega_2 
 - 
 \frac{\mathrm{i}}{2\pi}\sum_{j=1}^{N}Q_j
 \left\{
 \log(-q^{-1}) - 2\pi\mathrm{i}\frac{z-\zeta_j}{\omega_1} - u_j\omega_2
 \right\}
 \\
 = \: & 
 U\omega_2 + \sum_{j=1}^{N}Q_j\left(\frac{\zeta_j}{\omega_1}+\frac{\mathrm{i}}{2\pi}\omega_2 u_j\right), 
\end{align*}
where (\ref{eq:theta1-pseudo-periodicity}) is used on the second equality and (\ref{eq:q-sum-zero}) 
on the third equality. 
Then, the complex velocity $f_N^{\prime}(z) = u_N - \mathrm{i}v_N$ satisfies the double periodicity
\begin{equation}
 f_N^{\prime}(z+\omega_1) = f_N^{\prime}(z), \quad 
  f_N^{\prime}(z+\omega_2) = f_N^{\prime}(z),
\end{equation}
which means that $f_N^{\prime}(z)$ is an elliptic function of periods $\omega_1$ and $\omega_2$.

The approximate velocity potential $f_N(z)$ in (\ref{eq:mfs-solution}) 
is an analytic function in $\mathscr{D}$ such that it satisfies the condition 
(\ref{eq:uniform-flow2}) and gives the complex velocity $f_N^{\prime}(z)$ which is an elliptic function 
of periods $\omega_1$ and $\omega_2$. 
Regarding the boundary condition (\ref{eq:boundary-condition}), 
we pose a collocation condition on $f_N(z)$, namely, 
we assume the equation 
\begin{equation}
 \label{eq:collocation-cond}
  \im f_N(z_i) = C, \quad i = 1, \ldots, N, 
\end{equation}
where $z_1, \ldots, z_N$ are given points on $\partial D_{00}$ 
called the ^^ ^^ collocation points'' 
and $C$ is an unknown real constant. 
The equations (\ref{eq:collocation-cond}) are rewritten as 
\begin{multline}
 \label{eq:collocation-cond2}
  - \frac{1}{2\pi}\sum_{j=1}^{N}Q_j
  \left\{
  \log\left|\vartheta_1\left(\left.\frac{z_i-\zeta_j}{\omega_1}\right|\tau\right)\right|
  - 
  \re(u_j z_i)
  \right\}
  - C 
  = 
  - U(\im z_i), \\ i = 1, \ldots, N.
\end{multline}
The equations (\ref{eq:q-sum-zero}) and (\ref{eq:collocation-cond2}) form a system of 
linear equations with respect to the unknowns $Q_1, \ldots, Q_N$ and $C$. 
We determine the unknown charges $Q_j$ by solving this linear system 
and obtain the approximate velocity potential $f_N(z)$. 
Owing to the pseudo-periodicity (\ref{eq:pseudo-periodicity-potential1}) and 
(\ref{eq:pseudo-periodicity-potential2}), 
the approximate potential $f_N(z)$ automatically satisfies the collocation condition 
$\im f=\mbox{constant}$ on the boundary of other obstacle $D_{mn}$, $m,n\in\mathbb{Z}$, 
that is, 
\begin{equation*}
 \im f_N(z_i+m\omega_1+n\omega_2) = \mbox{constant}, \quad i = 1, \ldots, N.
\end{equation*}
\section{Numerical examples}
\label{sec:numerical-example}
In this section, we show some numerical examples which show the effectiveness of our method. 
All the computations were performed using programs coded in C++ with double precision. 

\begin{figure}[htbp]
 \begin{center}
  \begin{tabular}{cc}
   \psfrag{x}{\footnotesize$\re z/r$}
   \psfrag{y}{\footnotesize$\im z/r$}
   \includegraphics[width=0.49\textwidth]{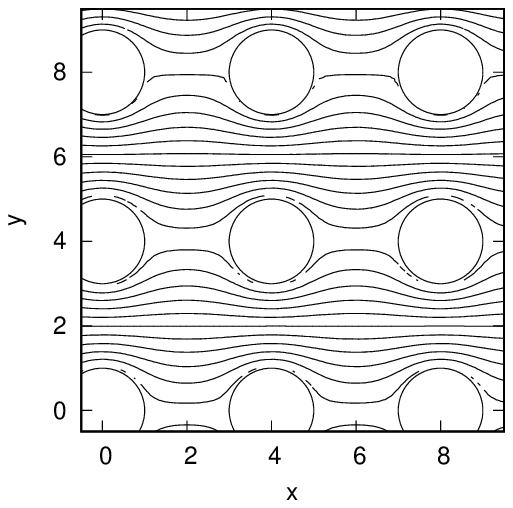} & 
       \psfrag{x}{\footnotesize$\re z/r$}
       \psfrag{y}{\footnotesize$\im z/r$}
       \includegraphics[width=0.49\textwidth]{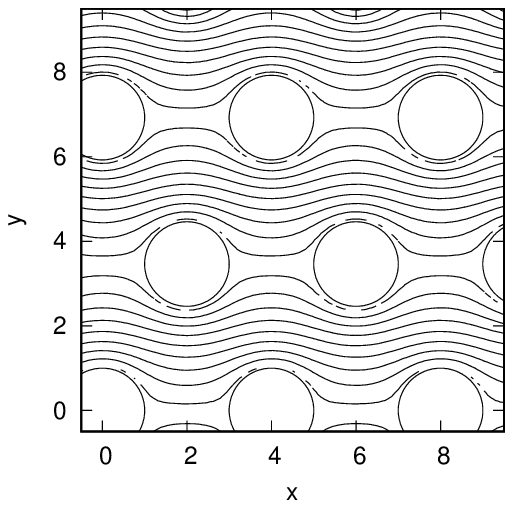}
       \\
   $(\omega_1,\omega_2)=(4r,4r\mathrm{i})$ & 
       $(\omega_1,\omega_2)=(4r,4r\mathrm{e}^{\mathrm{i}\pi/3})$
       \\
   \psfrag{x}{\footnotesize$\re z/r$}
   \psfrag{y}{\footnotesize$\im z/r$}
   \includegraphics[width=0.49\textwidth]{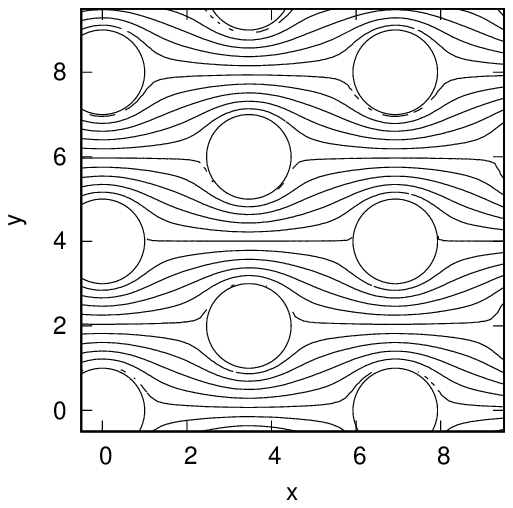}
   & 
       \psfrag{x}{\footnotesize$\re z/r$}
       \psfrag{y}{\footnotesize$\im z/r$}
       \includegraphics[width=0.49\textwidth]{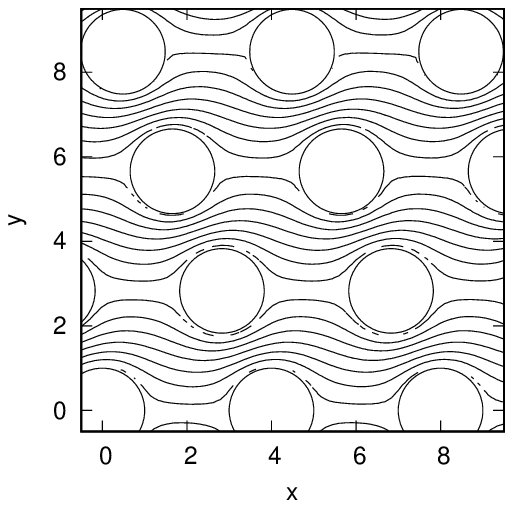}
   \\
   $(\omega_1,\omega_2)=(4r\mathrm{e}^{\mathrm{i}\pi/6},4r\mathrm{e}^{\mathrm{i}\pi/2})$ & 
       $(\omega_1,\omega_2)=(4r,4r\mathrm{e}^{\mathrm{i}\pi/4})$
  \end{tabular}
 \end{center}
 \caption{The streamlines of potential flows past a doubly-periodic array of cylinders with periods 
 $\omega_1, \omega_2$.}
 \label{fig:example1}
\end{figure}

Figure \ref{fig:example1} shows the examples of flows past a doubly-periodic array of cylinders, 
that is, flows in the domain
\begin{gather*}
 \mathscr{D} = 
 \mathbb{C} \setminus \bigcup_{m,n\in\mathbb{Z}}\overline{D_{mn}},
 \intertext{where}
 D_{mn} = 
 \left\{ \: z\in\mathbb{C} \: \left| \: |z-m\omega_1-n\omega_2| < r \: \right.\right\}, 
 \quad 
 m, n \in \mathbb{Z}
\end{gather*}
with a positive constant $r$ and periods $\omega_1, \omega_2\in\mathbb{C}$ such that 
$\im(\omega_2/\omega_1)>0$ and 
$\overline{D_{mn}}\cup\overline{D_{m^{\prime}n^{\prime}}}=\emptyset$  
if $(m,n)\neq(m^{\prime},n^{\prime})$. 
The charge points $\zeta_j$ and the collocation points $z_j$ are respectively taken as
\begin{equation}
 \label{eq:charge-collocation-point}
  \zeta_j = qr\exp\left(\mathrm{i}\frac{2\pi(j-1)}{N}\right), \quad 
  z_j = r\exp\left(\mathrm{i}\frac{2\pi(j-1)}{N}\right), 
  \quad j = 1, \ldots, N, 
\end{equation}
where $q$ is a constant such that $0 < q < 1$ and was taken as $q=0.7$ in the examples, 
and $u_j$ in (\ref{eq:uj}) are computed by the Monte Carlo method with a million points in $\mathscr{D}_0$.

\begin{figure}[htbp]
 \begin{center}
  \begin{tabular}{cc}
   \psfrag{n}{$N$}
   \psfrag{e}{$\epsilon_N$}
   \includegraphics[width=0.45\textwidth]{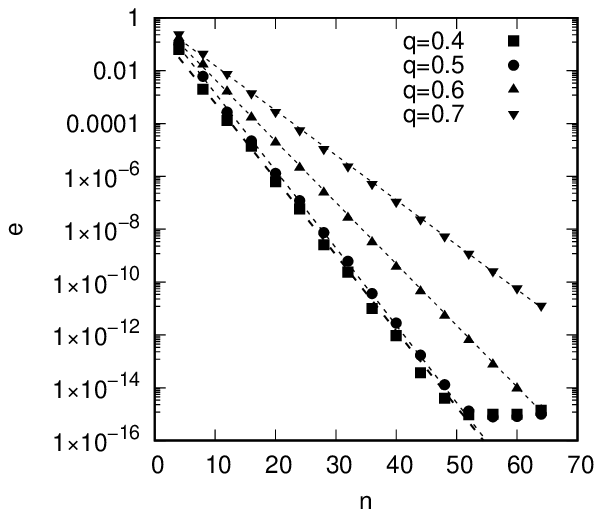}
   & 
       \psfrag{n}{$N$}
       \psfrag{e}{$\epsilon_N$}
       \includegraphics[width=0.45\textwidth]{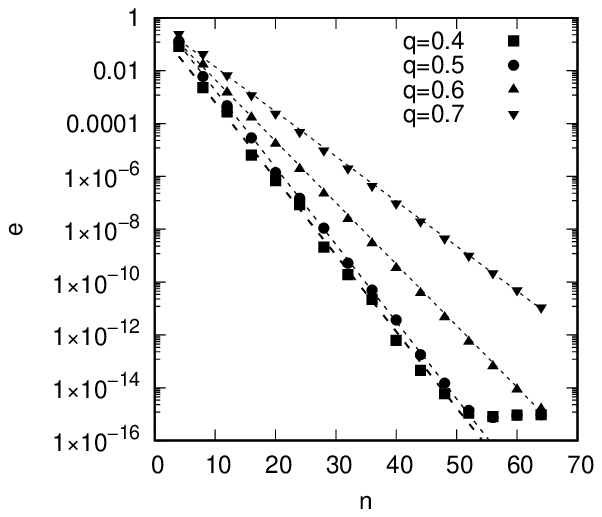}
       \\
   $(\omega_1,\omega_2)=(4r,4r\mathrm{i})$ & $(\omega_1,\omega_2)=(4r,4r\mathrm{e}^{\mathrm{i}\pi/3})$
       \\
   \psfrag{n}{$N$}
   \psfrag{e}{$\epsilon_N$}
   \includegraphics[width=0.45\textwidth]{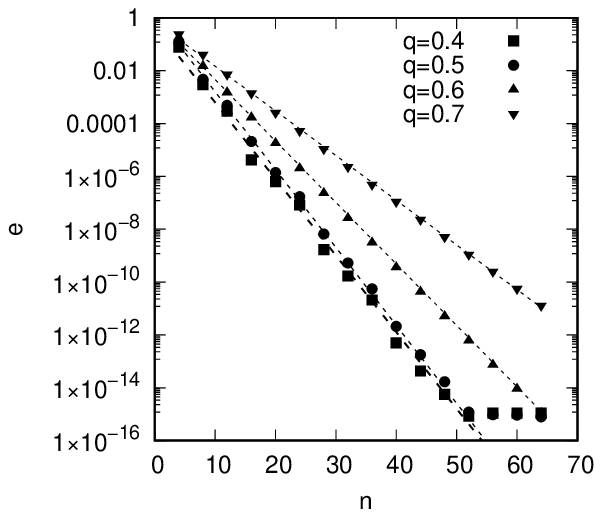}
   & 
       \psfrag{n}{$N$}
       \psfrag{e}{$\epsilon_N$}
       \includegraphics[width=0.45\textwidth]{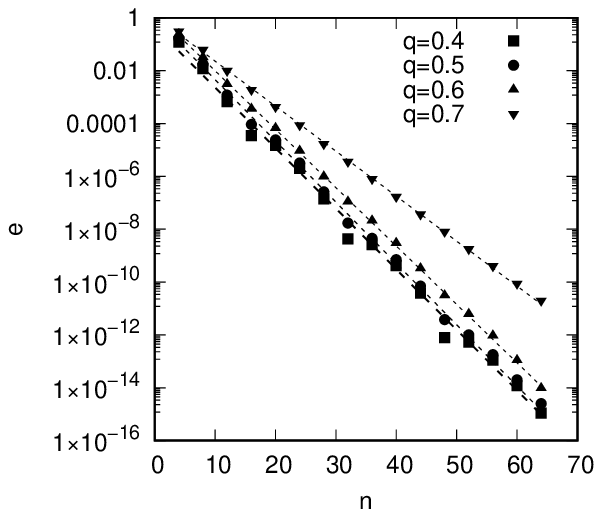}
       \\
   $(\omega_1,\omega_2)=(4r\mathrm{e}^{\mathrm{i}\pi/6})$ & 
       $(\omega_1,\omega_2) = (4r,4r\mathrm{e}^{\mathrm{i}\pi/4})$
  \end{tabular}
 \end{center}
 \caption{The error estimate $\epsilon_N$ of the method of fundamental solutions. }
 \label{fig:error}
\end{figure}

To estimate the accuracy of our method, we evaluated the value
\begin{equation*}
 \epsilon_N = 
  \frac{1}{Ur}\max_{z\in\partial D_{00}}|\im f_N(z) - C|, 
\end{equation*}
where $C$ is the constant obtained in solving the system of linear equations 
(\ref{eq:q-sum-zero}) and (\ref{eq:collocation-cond}). 
The value $\epsilon_N$ shows how accurately the approximate potential $f_N(z)$ satisfies the boundary condition 
(\ref{eq:boundary-cond0}). 
Figure \ref{fig:error} shows $\epsilon_N$ evaluated for the above numerical example with 
$\omega_1 = 4r$ and $\omega_2 = 4r\mathrm{i}$ for $q=0.4, 0.5, 0.6, 0.7$, 
where $q$ is the parameter appearing in (\ref{eq:charge-collocation-point}). 
The figure shows that the approximation of our method converges exponentially as the number of unknowns $N$ increases. 
Table \ref{tab:error} shows the decay rates of the error estimates $\epsilon_N$, 
which are also shown in Figure \ref{fig:error} in broken lines. 
The table shows that the error estimate $\epsilon_N$ roughly obeys 
\begin{equation*}
 \epsilon_N = \mathrm{O}(q^N)
\end{equation*}
for $q=0.5, 0.6, 0.7$ in the examples of 
$(\omega_1,\omega_2) = (4r, 4r\mathrm{i}), (4r,4r\mathrm{e}^{\mathrm{i}\pi/3}), 
(4r\mathrm{e}^{\mathrm{i}\pi/6}, 4r\mathrm{i})$ 
and for $q=0.6, 0.7$ in the example of $(\omega_1,\omega_2)=(4r,2r(1+\mathrm{i}))$. 

\begin{table}[htbp]
 \caption{The decay rates of the error estimates $\epsilon_N$.}
 \begin{center}
  \begin{tabular}{ccccc}
   \hline\rule{0pt}{12pt}
   $(\omega_1, \omega_2)$ & \multicolumn{4}{c}{$q$} 
       \\
   \cline{2-5}
   & \rule{0pt}{12pt}$0.4$ & $0.5$ & $0.6$ & $0.7$ 
   \\
   \hline\rule{0pt}{12pt}
   $(4r,4r\mathrm{i})$ & 
       $\mathrm{O}(0.52^N)$ & $\mathrm{O}(0.51^N)$ & $\mathrm{O}(0.59^N)$ & $\mathrm{O}(0.68^N)$ 
   \\ 
   $(4r,4r\mathrm{e}^{\mathrm{i}\pi/3})$ & 
       $\mathrm{O}(0.51^N)$ & $\mathrm{O}(0.51^N)$ & $\mathrm{O}(0.58^N)$ & $\mathrm{O}(0.68^N)$
   \\
   $(4r\mathrm{e}^{\mathrm{i}\pi/6}, 4r\mathrm{i})$ & 
       $\mathrm{O}(0.51^N)$ & $\mathrm{O}(0.51^N)$ & $\mathrm{O}(0.59^N)$ & $\mathrm{O}(0.68^N)$
   \\
   $(4r,2r(1+\mathrm{i}))$ & 
       $\mathrm{O}(0.59^N)$ & $\mathrm{O}(0.59^N)$ & $\mathrm{O}(0.60^N)$ & $\mathrm{O}(0.68^N)$
   \\
   \hline
  \end{tabular}
 \end{center}
 \label{tab:error}
\end{table}

In the actual computations, 
the condition (\ref{eq:uniform-flow}) that the average of the velocity is $\langle\boldsymbol{v}\rangle=(U,0)$ 
is not exactly satisfied because the integrals on $\mathscr{D}_0$ giving $u_j$ by (\ref{eq:mfs-solution}) 
are approximately computed by the Monte-Carlo method. 
We computed the average velocity $\langle\boldsymbol{v}\rangle$ by evaluating the integral in (\ref{eq:uniform-flow2}) 
using the Monte-Carlo method with a million points. 
Table \ref{tab:average-velocity} shows the result. 
The table shows that the condition (\ref{eq:uniform-flow}) is approximately satisfied 
with error $\sim 10^{-4}$ to $10^{-5}$. 
\begin{table}[htbp]
 \caption{The actual average velocities $\langle\boldsymbol{v}\rangle$.}
 \begin{center}
  \begin{tabular}{cc}
   \hline \rule{0pt}{12pt}
   $(\omega_1,\omega_2)$ & $\langle\boldsymbol{v}\rangle/U$ \\
   \hline \rule{0pt}{12pt}
   $(4r,4r\mathrm{i})$ & $(1.0004, -9\times 10^{-5})$ \\ 
   $(4r,4r\mathrm{e}^{\mathrm{i}\pi/3})$ & $(0.9997,-2\times 10^{-4})$ \\ 
   $(4r\mathrm{e}^{\mathrm{i}\pi/6}, 4r\mathrm{i})$ & $(0.9998,3\times 10^{-5})$\\ 
       $(4r,2r(1+\mathrm{i}))$ & $(0.9998, 3\times 10^{-5})$
       \\
   \hline
  \end{tabular}
 \end{center}
 \label{tab:average-velocity}
\end{table}
\section{Concluding Remarks}
\label{sec:concluding-remark}
In this paper, we proposed a method of fundamental solutions for 
the problems of two-dimensional potential flow past a doubly-periodic array of obstacles. 
In terms of mathematics, our problem is to find the complex velocity potential, 
an analytic function in the flow domain with double periodicity. 
The solution obviously involves a doubly-periodic function, and it is difficult to approximate it 
by the conventional method. 
Then, we proposed a new method of fundamental solution for this problem using the periodic 
fundamental solutions, which is the logarithmic potential with a doubly-periodic array of sources 
constructed using the theta functions. 
The proposed method inherits the advantages of the conventional method of fundamental solutions 
and approximates the solution of our problem with double periodicity well. 
The numerical examples showed the effectiveness of our method. 

We have two issues regarding this paper for future study. 
The first problem is a theoretical error estimate of our method. 
Theoretical studies on the accuracy of the method of fundamental solution have been presented for 
special cases such as two-dimensional Laplace equation 
and Helmholtz equation in a disk \cite{KatsuradaOkamoto1988,ChibaUshijima2009,OgataChibaUshijima2011} 
and two-dimensional Laplace equation in a domain with an analytic boundary 
\cite{Katsurada1990,OgataKatsurada2014}. 
However, theoretical error estimate still remains unknown for many types of problems and methods including 
our method. 
The author believes it one of the most important works on the method of fundamental solutions 
to give a theoretical error estimate.

The second problem is to extend our method to other problems with periodicity than two-dimensional potential problems. 
The author has already given methods of fundamental solutions for Stokes flow problems 
\cite{OgataAmano2010,OgataAmanoSugiharaOkano2003} and elasticity problems \cite{Ogata2008} with periodicity. 
In the work on periodic Stokes flow \cite{OgataAmanoSugiharaOkano2003}, 
the author et al. proposed to use the periodic fundamental solutions, which was presented by Hasimoto 
\cite{Hasimoto1959} and is given by a Fourier series. 
It is expected to construct a method of fundamental solutions 
using periodic fundamental solutions given by the theta functions 
or elliptic functions \cite{Hasimoto2009,CrowdyLuca2018} as in this paper. 
\section*{Acknowledgement}
 The author thanks Professor Yusaku Yamamoto for his advice in numerical computations. 
\bibliographystyle{plain}
\bibliography{arxiv2020_01}

\begin{thebibliography}{10}

\bibitem{Amano1994}
K.~Amano.
\newblock A charge simulation method for the numerical conformal mapping of
  interior, exterior and doubly-connected domains.
\newblock {\em J. Comput. Appl. Math.}, 53(3):353--370, 1994.

\bibitem{Amano1998}
K.~Amano.
\newblock A charge simulation method for numerical conformal mapping onto
  circular and radial slit domains.
\newblock {\em SIAM J. Sci. Comput.}, 19(4):1169--1187, 1998.

\bibitem{AmanoOkanoOgataSugihara2012}
K.~Amano, D.~Okano, H.~Ogata, and M.~Sugihara.
\newblock Numerical conformal mapping onto the linear slit domains.
\newblock {\em Japan J. Indust. Appl. Math.}, 29:165--186, 2012.

\bibitem{Armitage-Eberlein2006}
J.~V. Armitage and W.~F. Eberlein.
\newblock {\em Elliptic Functions}.
\newblock Cambridge University Press, Cambridge, 2006.

\bibitem{ChibaUshijima2009}
F.~Chiba and T.~Ushijima.
\newblock Exponential decay of errors of a fundamental solution method applied
  to a reduced wave problem in the exterior region of a disc.
\newblock {\em J. Comput. Appl. Math.}, 231:869--885, 2009.

\bibitem{CrowdyLuca2018}
D.~Crowdy and E.~Luca.
\newblock Fast evaluation of the fundamental singularities of two-dimensional
  doubly periodic {Stokes} flow.
\newblock {\em J. Eng. Math.}, 111:95--110, 2018.

\bibitem{FairweatherKarageorghis1998}
G.~Fairweather and A.~Karageorghis.
\newblock The method of fundamental solutions for elliptic boundary value
  problems.
\newblock {\em Adv. Comp. Math.}, 9:69--95, 1998.

\bibitem{GreengardKropinski2004}
L.~Greengard and M.~C. Kropinski.
\newblock Integral equation methods for stokes flow in doubly-periodic domains.
\newblock {\em J. Eng. Math.}, 48:157--170, 2004.

\bibitem{Hasimoto1959}
H.~Hasimoto.
\newblock On the periodic fundamental solutions of the stokes equations and
  their application to viscous flow past a cubic array of spheres.
\newblock {\em J. Fluid. Mech.}, 5(2):317--328, 1959.

\bibitem{Hasimoto2008}
H.~Hasimoto.
\newblock Periodic fundamental solution of a two-dimensional {Poisson}
  equation.
\newblock {\em J. Phys. Soc. Japan}, 77(10):104601, 2008.

\bibitem{Hasimoto2009}
H.~Hasimoto.
\newblock Periodic fundamental solution of the two-dimensional {Stokes}
  equations.
\newblock {\em J. Phys. Soc. Japan}, 78(7):074401, 2009.

\bibitem{Katsurada1990}
M.~Katsurada.
\newblock Asymptotic error analysis of the charge simulation method in a jordan
  region with an analytic boundary.
\newblock {\em J. Fac. Sci. Univ. Tokyo, Sect. IA Math.}, 37:635--657, 1990.

\bibitem{KatsuradaOkamoto1988}
M.~Katsurada and H.~Okamoto.
\newblock A mathematical study of the charge simulation method {I}.
\newblock {\em J. Fac. Sci. Univ. Tokyo, Sect IA}, 35(3):507--518, 1988.

\bibitem{Liron1978}
N.~Liron.
\newblock Fluid transport by cilia between parallel plates.
\newblock {\em J. Fluid Mech.}, 86(4):705--726, 1978.

\bibitem{Milne-Thomson2011}
L.~M. Milne-Thomson.
\newblock {\em Theoretical Hydrodynamics}.
\newblock Dover, New York, 2011.

\bibitem{Murashima1983}
S.~Murashima.
\newblock {\em Charge Simulation Method and Its Applications}.
\newblock Morikita-Shuppan, Tokyo, 1983.
\newblock (in Japanese).

\bibitem{Murota1993}
K.~Murota.
\newblock On {^^ ^^ invariance"} of schemes in the fundamental solution method.
\newblock {\em Trans. IPS Japan}, 34(3):533--535, 1993.
\newblock (in Japanese).

\bibitem{Murota1995}
K.~Murota.
\newblock Comparison of conventional and {^^ ^^ invariant"} schemes of
  fundamental solutions method for annular domains.
\newblock {\em Japan J. Indust. Appl. Math.}, 12:61--85, 1995.

\bibitem{Ogata2008}
H.~Ogata.
\newblock Fundamental solution method for periodic plane elasticity.
\newblock {\em J. Numer. Anal. Indust. Appl. Math. (JNAIAM)}, 3(3--4):249--267,
  2008.

\bibitem{OgataAmano2010}
H.~Ogata and K.~Amano.
\newblock Fundamental solution method for two-dimensional stokes flow problems
  with one-dimensional periodicity.
\newblock {\em Japan J. Indust. Appl. Math.}, 27:191--215, 2010.

\bibitem{OgataAmanoSugiharaOkano2003}
H.~Ogata, K.~Amano, M.~Sugihara, and D.~Okano.
\newblock A fundamental solution method for viscous flow problems with
  obstacles in a periodic array.
\newblock {\em J. Comput. Appl. Math.}, 152(1--2):411--425, 2003.

\bibitem{OgataChibaUshijima2011}
H.~Ogata, F.~Chiba, and T.~Ushijima.
\newblock A new theoretical error estimate of the method of fundamental
  solutions applied to reduced wave problems in the exterior region of a disk.
\newblock {\em J. Comput. Appl. Math.}, 235(12):3395--3412, 2011.

\bibitem{OgataKatsurada2014}
H.~Ogata and M.~Katsurada.
\newblock Convergence of the invariant scheme of the method of fundamental
  solutions for two-dimensional potential problems in a jordan region.
\newblock {\em Japan J. Indust. Appl. Math.}, 31:231--262, 2014.

\bibitem{OgataOkanoAmano2002}
H.~Ogata, D.~Okano, and K.~Amano.
\newblock Numerical conformal mapping of periodic structure domains.
\newblock {\em Japan J. Indust. Appl. Math.}, 19:257--275, 2002.

\bibitem{Sanchez-SezmaRosenblueth1979}
F.~J. Sanchez-Sezma and E.~Rosenblueth.
\newblock Ground motion at canyons of arbitrary shape under incident sh waves.
\newblock {\em Int. J. Earthq. Eng. Struct. Dyn.}, 7:441--450, 1979.

\bibitem{SingerSteinbiglerWeiss1974}
H.~Singer, H.~Steinbigler, and P.~Weiss.
\newblock A charge simulation method for the calculation of high voltage
  fields.
\newblock {\em IEEE Trans. Power Appar. Syst.}, PAS-93:1660--1668, 1974.

\bibitem{Steinbigler-dissertation1969}
H.~Steinbigler, 1969.
\newblock dissertation, Tech. Univ. {M\"{u}nchen}.

\bibitem{ZickHomsy1982}
A.~A. Zick and G.~M. Homsy.
\newblock Stokes flow through periodic arrays of spheres.
\newblock {\em J. Fluid Mech.}, 115:13--26, 1982.

\end{thebibliography}
\end{document}